\documentclass[11pt]{article}
\usepackage{graphicx, color, amssymb}
\def\ds{\displaystyle}
\def\forall{\hbox{for all}~}

\def\H{{\cal{H}}}

\def\S{{\cal S}}

\def\I{{\cal I}}

\def\P{{\cal P}}
\def\E{{\cal E}}

\def\ve{\varepsilon}

\def\R{\mathbb R}

\def\wto{\rightharpoonup}

\def\c{\centerline}
\def\P{{\cal P}}

\def\vs{\vskip 2em}

\def\v{\vskip 1em}
\def\O{{\cal O}}

\def\begi{\begin{itemize}}
\def\endi{\end{itemize}}

\def\ov{\overline}
\def\Tilde{\widetilde}
\def\Hat{\widehat}
\def\bega{\begin{array}}
\def\enda{\end{array}}
\def\meas{\hbox{meas}}

\def\bel{\begin{equation}\label}
\def\eeq{\end{equation}}
\def\sqr#1#2{\vbox{\hrule height .#2pt
\hbox{\vrule width .#2pt height #1pt \kern #1pt
\vrule width .#2pt}\hrule height .#2pt }}
\def\square{\sqr74}
\def\endproof{\hphantom{MM}\hfill\llap{$\square$}\goodbreak}

\newtheorem{theorem}{Theorem}[section]
\newtheorem{lemma}[theorem]{Lemma}

\newtheorem{definition}[theorem]{Definition}
\newtheorem{remark}[theorem]{Remark}

\begin{document}

\title{\bf  Irrigable Measures for Weighted Irrigation Plans}

\author{ Qing Sun
\\
\,
\\
Department of Mathematics, Penn State University \\
e-mails: qxs15@psu.edu
}

\maketitle
\begin{abstract} 
	A model of irrigation network, where lower branches must be thicker in order to support the weight of the higher ones, was recently introduced in \cite{BSUN2}.  This leads to a countable family of ODEs, describing the thickness of every branch, solved by backward induction.  The present paper determines what kind of measures can be irrigated with a finite weighted cost.  Indeed, the boundedness of the cost depends on  the dimension of the support of the irrigated measure, and also on the asymptotic properties of the ODE which determines the thickness of branches.
	
\end{abstract}
\vs

\section{Introduction}
\label{s:1}
\setcounter{equation}{0}

In a ramified transport network \cite{BCM,BCM1,MMS,X3,X15,MS},  the 
Gilbert transport cost along each arc is computed by 
\bel{1} [\hbox{length}] \times [\hbox{flux}]^\alpha   \eeq
for some given $\alpha\in [0, 1]$. 
When $\alpha< 1$, this accounts for an economy of scale:  transporting the same amount of particles is cheaper if these particles travel together along the same arc.

In the recent paper \cite{BSUN2}, the authors considered an irrigation plan where the cost
per unit length is determined by a weight function $W$ .
The main motivation behind this model is that, for a free standing structure like a tree, the lower portion of each branch needs to bear the weight of the upper part.
Hence, even if the flux of water and nutrients is constant along a branch, the thickness
(and hence the cost per unit length) grows as one moves from the tip toward the root. In the variational problems of optimal tree roots and branches\cite{BPSUN,BSUN}, this ``weighted irrigation cost" is more suitable to model the associated cost for transporting water and nutrients from the roots to the leaves.

In this model, the weights are constructed inductively, starting from the outermost  branches
and proceeding toward the root. Along each branch, the weight $W$ is determined by
solving a suitable
 ODE, possibly with measure-valued right hand side.  This is more conveniently written
in the integral form
\bel{qs0}W(s)~=~\int_s^\ell f(W(\sigma))\, d\sigma + m(s),\eeq
where $s\in [0,\ell]$ is the arc-length parameter along the branch, $s\mapsto m(s)$ 
is a non-increasing function describing the flux, and $f$ is a non-negative, continuous function. A natural set of assumptions on $f$ is 
\begin{itemize}

\item[{\bf (A1)}] {\it The function $f: \R_+\mapsto \R_+$ is continuous on $[0,+\infty[\, $, twice continuously differentiable for $z>0$, and satisfies
	\bel{qs00}f(0)~ = ~0,\quad f'(z)~>~0,\quad f''(z)~\leq~ 0\qquad \forall ~z >0 .\eeq}
	
\end{itemize}

The main result in \cite{BSUN2}
established the lower semicontinuity
 of the weighted irrigation cost, w.r.t.~the pointwise convergence of irrigation plans. 
In particular, for any positive, bounded Radon measure $\mu$, if there is an 
 admissible irrigation plan whose weighted cost is finite,
 then there exists an irrigation plan for $\mu$ with minimum cost.

The goal of the present paper is to understand  whether 
 a given Radon measure $\mu$ irrigable or not, with respect to the weighted irrigation cost. 
 That is, whether there exists an irrigation plan for $\mu$ whose weighted irrigation cost is finite.  In the case without weights, i.e., with the classical Gilbert cost (\ref{1}), this problem has 
 been studied in \cite{DS}, and further investigated in \cite{A,B,D}.  The authors in \cite{DS} proved that if a measure $\mu$ is $\alpha$-irrigable, 
 then it  must be concentrated on a set with Hausdorff dimension  $\leq {1\over 1 -\alpha}$. 
On the other hand,  if $\alpha> 1 - {1\over d}$, 
 then every bounded Radon measure with bounded support in $\R^d$ 
 has finite irrigation cost \cite{BCM,DS}. 
 
 As shown by our analysis, in the presence of weights the irrigability of a measure $\mu$
 depends on the dimension of the set where $\mu$ is concentrated, on
 the exponent $\alpha$,  and also on the asymptotic behavior of the function $f(z)$ as $z\to 0^+$.   

The remainder of the paper is organized as follows. Section 2 reviews the construction of the weight functions on the various branches of an  irrigation plan. In Section 3 we prove our main results on the irrigability of Radon measures. 

\section{Review of the weighted irrigation plans}\label{s:review}
\setcounter{equation}{0}

\subsection{Weight functions on finitely many branches}\label{finite}

\begin{figure}[ht]
	\centering
	\includegraphics[width=10cm]{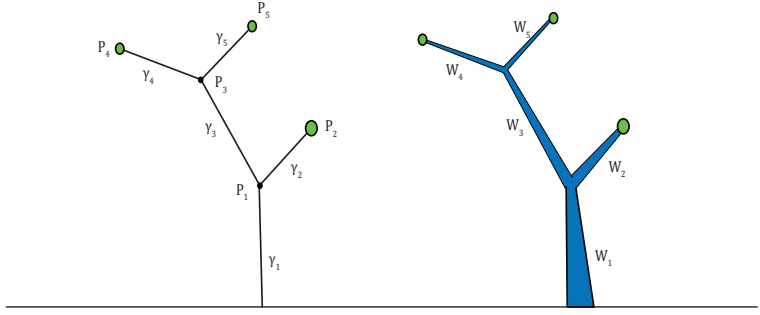}
	\caption{\small Left: A free standing tree with 5 branches. In this example, $\O(1) = \{2,3 \},\O(3) = \{4,5\}, \O(2)=\O(4)=\O(5)=\emptyset$. Right: On each branch, the weight decreases as one moves from the lower portion to the tip.}
	\label{f:e1}
\end{figure}

To illustrate the basic idea of the weighted irrigation model, we first consider a  network with finitely many branches. 
As shown on the left of Fig. \ref{f:e1}, each directed branch will be denoted by 
$\gamma_i:[a_i,b_i]\mapsto \R^d, i = 1,\ldots,N$, oriented from the root toward the tip and parameterized by arc-length. Call $P_i = \gamma_i(b_i)$ the ending node of the branch $\gamma_i$. 

On each branch $\gamma_i$, we first prescribe a  left-continuous, non-increasing function $m_i:[a_i,b_i]\mapsto \R^+$, which can be interpreted as the ``flux" along the branch. 
Roughly speaking, $m_i(t)$ is the amount of mass transported through the point $\gamma_i(t)$.

 Call
$\O(i)$
the set of index labelling the branches that originate from the node $P_i=\gamma_i(b_i)$, that is
\bel{12}\O(i)~=~\Big\{ j\in \{1,\ldots, N\};~\gamma_j(a_j) = P_i \Big\}.\eeq
Moreover, consider the sets of indices inductively defined by
\begin{equation}\label{qs5}
\bega{c}
\ds  \I_1 ~\doteq~\left\{i\in\{1,\ldots,N\}\, ;~\O(i) = \emptyset \right\},\\[3mm]
\ds \I_{k+1}~\doteq~\left\{i \in \{1,\ldots,N\}\, ;~\O(i) \subseteq \I_1\cup\cdots\cup \I_k \right\}\setminus (\I_1\cup\cdots\cup \I_k).
 \enda
\end{equation}

From \cite{BSUN2} the weight function $W_i(\cdot)$ on each branch $\gamma_i$ is defined inductively on $\I_k, k\geq 1$.
\begin{itemize}
	\item[(i)] For $k=1$, on each branch $\gamma_i:[a_i,b_i]\mapsto \R^d$ with $i\in \I_1$, the weight $W_i:[a_i,b_i]\mapsto \R_+$ is defined to be the solution of 
	\bel{qs1}\omega(t)~=~\int_t^{b_i} f(\omega(s)) \, ds + m_i(t),\qquad t\in\, ]a_i, b_i],\eeq
	where $f$ is a given function,  satisfying {\bf (A1)}, and $m_i$ is the flux along the branch.
	
	\item[(ii)] Assume the weight functions $W_i(t)$ have already been constructed along all branches $\gamma_j:[a_j,b_j]\mapsto \R^d$ with $j\in \I_1\cup\ldots\cup \I_{k-1}$.
	
	For $i\in \I_k$, the weight $W_i(t)$ along the $i$-th branch is defined to be the solution of 
	\bel{qs2}
	\omega(t)~=~\int_t^{b_i} f(\omega(s))\, ds + m_i(t) + \ov\omega_i,\qquad t\in \, ]a_i, b_i],
	\eeq
	where 
	\bel{qs3} 
	\ov\omega_i~\doteq~\sum_{j\in \O(i)} W_j(a_j^+) - \sum_{j\in \O(i)} m_j(a_j^+).
	\eeq
	
\end{itemize}


\subsection{Irrigation plans for general measures}
Following Maddalena, Morel, and Solimini \cite{MMS}, the transport network for general Radon measure can be described in  a Lagrangian way. 
 Let $\mu$ be a fixed Radon measure on $\R^d$ with $\mu(\R^d) = M$ and let $\Theta = [0,M]$. We think of 
 $\theta\in \Theta$ as a Lagrangian variable, labelling a water particle.
  An {\bf irrigation plan for $\mu$} is a function 
$$\chi: \Theta\times \R_+ \mapsto \R^d,$$
measurable w.r.t.~$\theta$ and continuous w.r.t.~$t$, which satisfies the following conditions:
\begin{itemize}
	
		\item All particles initially lie at the origin: $\chi(\theta,0) = 0~ \forall \theta\in \Theta$.
		
	\item  For a.e.~$\theta\in \Theta$ the map $t\mapsto  \chi(\theta,t)$ 
	is 1-Lipschitz  and constant for $t$ large.  
	Namely, there exists $\tau(\theta)\geq 0$ such that
	$$\left\{\bega{cl} |\chi(\theta,t)-\chi(\theta, s)|~\leq~|t-s|\qquad&\hbox{for all}~ t,s\geq 0,\\[3mm]
	\chi(\theta, t)~=~\chi(\theta,\tau(\theta))\qquad &\hbox{for every} ~t\geq \tau(\theta).\enda\right.$$
	Throughout the following, $\tau(\theta)$ will denote the smallest time $\tau$ such that $\chi(\theta,\cdot)$ is constant for $t\geq \tau$.

			\item $\chi$ irrigates the measure $\mu$. That is, for each Borel set $V\subseteq \R^d$,
		$$ \mu(V)~=~\meas\left(\{ \theta\in \Theta;~\chi(\theta,\tau(\theta)) \in V  \}\right).  $$

\end{itemize}
One can think of $\chi(\theta,t)$ as  the position of particle $\theta$ at time $t$. 

To define the flux on $\chi$, which measures the total amount of particles travelling along the same path, we first need  an equivalence relation between two Lipschitz maps.

\begin{definition}  
	\label{d:eq1}
	We say that two 1-Lipschitz maps $\gamma:[0,t]\mapsto \R^d$ and $\Tilde\gamma:[0, \Tilde{t}\, ]\mapsto\R^d$ 
	are {\bf equivalent} if they are parametrizations of the same curve, and write it as $\gamma\simeq\Tilde\gamma$. When we use the arc-length re-parametrization 
	$$ \sigma\mapsto \gamma(s(\sigma)),\qquad\qquad \hbox{\rm where}\qquad  \int_0^{s(\sigma)} |\dot{\gamma}(t)|\, dt~=~\sigma,  $$
	then two 1-Lipschitz maps are equivalent means their arc-length re-parametrizations coincide.
\end{definition}

Throughout the following, 
we denote by $\gamma\Big|_{[0,t]}$ the restriction of a map $\gamma$ to the interval
$[0,t]$.

\begin{definition}\label{d:1}
	Let $\chi:\Theta\times \R_+\mapsto \R^d$ be an irrigation plan for the measure $\mu$. On the set $\Theta\times \R_+$, we write
 $(\theta,t)\,\sim\,(\theta', t')$
	whenever $\chi(\theta,\cdot)\Big|_{[0,t]}\simeq\chi(\theta',\cdot)\Big|_{[0,t']}$. This means that
	the maps
	$$s~\mapsto~\chi(\theta,s), \quad s\in [0,t]\qquad\hbox{and}
	\qquad s~\mapsto~\chi(\theta',s), \quad s\in [0,t']
	$$
	are equivalent in the sense of Definition~\ref{d:eq1}.
	
	The {\bf multiplicity} at $(\theta,t)$ is then defined as
	\bel{mtx}m(\theta,t)~\doteq~\meas\Big(\bigl\{ \theta'\in\Theta\,;~~(\theta', t')\sim (\theta,t) ~~\hbox{for some}~t'>0\bigr\}
	\Big).\eeq
\end{definition}

Given an irrigation plan $\chi:\Theta\times\R_+\mapsto\R^d$, in order to have finite weighted irrigation cost constructed in the next section, we should  always assume the following conditon.
\begi
\item[{\bf (A2)}] {\it
	For a.e.~$\theta\in \Theta$,  one has $m(\theta,t)>0$
	for every $0\leq t<\tau(\theta)$. }
\endi
\subsection{Weight functions for an irrigation plan.}

Given a bounded Radon measure $\mu$ in $\R^d$ and an irrigation plan  $\chi: \Theta\times \R_+\mapsto \R^d$ for $\mu$,  in this section we review the construction of the weight function $W=W(\theta,t)$ on the irrigation plan. Notice that for an irrigation plan $\chi$ of a general Radon measure,  for each particle $\theta \in \Theta$, the map $\chi(\theta,\cdot):\R_+\mapsto \R^d$ describes a continuous curve in $\R^d$. Thus $\chi$ may contain infinitely many branches. To construct the weight function on each branch, the idea is to first compute the weights $W^\ve$ on $\chi^\ve$, which is the truncation of $\chi$ on the branches with multiplicity $\geq \ve$. It turns out that $\chi^\ve$ only consists of finitely many branches, so that we can compute $W^\ve$ as in Section \ref{finite} . The weight $W$ is then constructed by taking the limit of $W^\ve$, as  $\ve\to 0+$.
\begin{definition}\label{d:good}
	Given an irrigation plan $\chi$, a  path $\gamma:[0, \ell]\mapsto\R^d$, parameterized by 
	arc-length,
	is {\bf $\ve$-good} if and only if
	\bel{goode}\meas\left(\Big\{\theta\in \Theta\,;~~\chi(\theta,\cdot)\Big|_{[0,t]} \simeq
	~\gamma ~~\hbox{for some}~ t= t(\theta)>0\Big\}
	\right)~\geq~\ve,
	\eeq
		where the equivalence relation $\simeq$ is given in Definition \ref{d:eq1}.

\end{definition}
In other words,  $\gamma$ is $\ve$-good if there is an amount  $\geq \ve$
of particles whose trajectory contains $\gamma$ as  initial 
portion.

\v
For any given $\ve>0$, following \cite{BSUN2} we define the $\ve$-stopping time $\tau_{\ve}:\Theta\mapsto \R_+$
by setting
\bel{5} \tau_{\ve}(\theta)~\doteq~\max \left\{ t\geq 0;~m(\theta,t) \geq \ve \right\}.  \eeq
Define the $\ve$-truncation $\chi^\ve$ of irrigation plan $\chi$ as
\bel{6}  \chi^\ve(\theta,t)~\doteq~\left\{\bega{rl} \chi(\theta,t)&\qquad \hbox{ if }~~ t< \tau_{\ve}(\theta)\\[3mm]
\chi(\theta,\tau_{\ve}(\theta))&\qquad \hbox{ if }~~ t \geq \tau_{\ve}(\theta) \enda\right.  \eeq
In other words, in the $\ve$-truncation $\chi^\ve$, only those paths in $\chi$ with multiplicity $\geq \ve$ are kept.  For any $\theta\in \Theta$, if $\tau_{\ve}(\theta)>0$, the $\ve$-good portion $\chi(\theta,\cdot)\Big|_{[0,\tau_\ve(\theta)]}$ of the path $t\mapsto \chi(\theta,\cdot)$ is included in $\chi_\ve$.  
 
Notice that the family of all curves parameterized by arc-length comes with a natural partial order.
Namely, given two maps $\gamma:[0,\ell]\mapsto\R^d$,
$\Tilde\gamma:[0,\Tilde\ell]\mapsto\R^d$,  
we write $\gamma\preceq \Tilde\gamma$  if $\ell\leq \Tilde\ell$ and $\gamma(s)=\Tilde\gamma(s)$ for
all $s\in [0,\ell]$.   In the family of all $\ve$-good paths in the irrigation plan $\chi$, we can thus find the maximal $\ve$-good paths, w.r.t the above partial order. 
As shown  in \cite{BSUN2}, the total number of maximal $\ve$-good paths in the irrigation plan $\chi$ is bounded by ${M\over \ve}$, where $M$ is the total mass of $\mu$. Therefore, the $\ve$-truncation $\chi^\ve$ is a network with finitely many branches, 
consisting of all  maximal  $\ve$-good paths in $\chi$.

For a fixed $\ve>0$, to compute the weight functions on the $\ve$-truncation $\chi^\ve$, we now let  $\{\Hat \gamma_1,\ldots,\Hat \gamma_\nu\}$ be the set of all 
maximal $\ve$-good paths.
Along each path $\Hat \gamma_i:[0,\hat \ell_i]\mapsto\R^d$ we define the {\bf multiplicity} $\Hat m_i: [0, \hat \ell_i]
\mapsto \R_+$
by setting 
\bel{qs4}\Hat m_i(s)~=~\meas\bigg( \left\{ \theta\in \Theta\,;~~
\hbox{there exists $t\geq 0$ such that }~\chi(\theta,\cdot)\Big|_{[0,t]}~\simeq~
\Hat \gamma_i\Big|_{[0,s]}\right\}\bigg).\eeq

Since two maximal paths may coincide on the initial portion and  bifurcate later,  we consider the bifurcation times
\bel{tauij}
\tau_{ij}~=~\tau_{ji}~\doteq~\max\,\Big\{ t\geq 0\,;~~\Hat\gamma_i(s)=\Hat\gamma_j(s)~~\forall
s\in [0,t]\Big\}.\eeq
For each maximal path $\Hat\gamma_i$, we split it into several elementary branches $\gamma_k$, by the  following Path Splitting Algorithm{\bf (PSA)}, which is first introduced in \cite{BSUN2}. 
\begi
\item[{\bf (PSA)}] For each $i\in\{1,\ldots,\nu\}$,  consider the set 
$$\{\tau_{i1}, \ldots,\tau_{i\nu}\} ~=~\{t_{i,1},\ldots, t_{i,N(i)}\},$$
where the times
\bel{tij}0\,=t_{i,0}\,< \,t_{i,1}\,<\,t_{i,2}\,<\,\cdots\,<\,t_{i,N(i)}\, <\,\hat \ell_i\eeq
provide an increasing arrangement of the set of  times $\tau_{ij}$ where the path
$\Hat\gamma_i$ splits apart from other maximal paths.
For each $k=1,\ldots, N(i)$, let $\gamma_{i,k}$
be the restriction of the maximal path $\Hat \gamma_i$ to the subinterval $[t_{i,k-1}, t_{i,k}]$.
The multiplicity function $m_{i,k}$ along this path is defined simply as
\bel{multj}
m_{i,k}(t)~=~\Hat m_i(t)\qquad\qquad t\in [t_{i,k-1}, t_{i,k}].\eeq

If $\tau_{ij}>0$, i.e.~if the two maximal paths $\Hat \gamma_i$ and 
$\Hat \gamma_j$ partially overlap, it is clear that  some of the elementary branches $\gamma_{i,k}$
will coincide with some $\gamma_{j,l}$.
To avoid listing multiple times the same branch, we thus remove from our list 
all branches
$\gamma_{j,l}:[t_{j, l-1}, t_{j,l}]\mapsto\R^d$ such that 
$t_{j,l} \leq \tau_{ij}$ for some $i<j$.    
After relabeling all the remaining branches, the algorithm yields a 
family of elementary branches and corresponding multiplicities 
\bel{eph}\gamma_i:[a_i,b_i]\mapsto \R^d,\qquad m_i:\, [a_i,b_i]\mapsto \R_+\,,
\qquad i=1,\ldots,N \eeq
where $N$ is the total number of elementary branches.
\endi

\begin{figure}[ht]
	\centering
	\includegraphics[width=15cm]{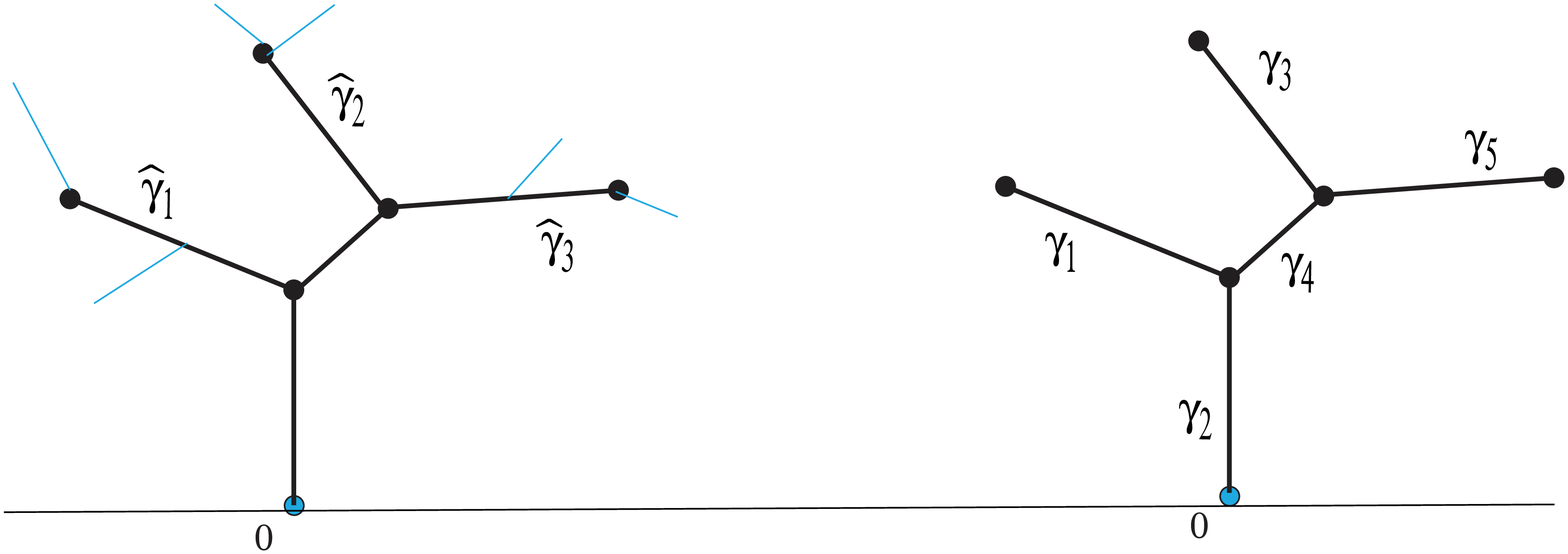}
	\caption{\small Left: Two finite truncation plans, showing three maximal $\ve$-good paths (thick lines) and six maximal $\ve'$-good paths (thin lines), for $0<\ve'<\ve$. Right: The three maximal $\ve$-good paths can be partitioned into five elementary branches, by the Path Splitting Algorithm. }
	\label{f:1}
\end{figure}

On these elementary branches $\gamma_i, i\geq 1$, we can compute the weight function $W_i$ on each $\gamma_i$ inductively, as in Section \ref{finite}.


On each maximal $\ve$-good path $\Hat\gamma_i$ with $1\leq i\leq \nu$,
the above construction yields a weight $\Hat W_{i,k}$ on the restriction 
of  $\Hat \gamma_i$ to each subinterval $[t_{i, k-1}, t_{i,k}]$.
Along the maximal path 
$\Hat\gamma_i$, the weight $\Hat W_i:[0, \hat\ell_i]\mapsto \R_+$
is  then defined simply by setting
\bel{HWI}\Hat W_i(t)~=~\Hat W_{i,k}(t)\qquad\qquad \hbox{if}~~~
t\in [t_{i, k-1}, t_{i,k}].\eeq

\v
Next, on the $\ve$-truncation $\chi^\ve$
we define the
weight function $W^\ve:\Theta\times
\R_+\mapsto \R_+$ by setting
\bel{We1}
W^\ve(\theta,t)~\doteq~\left\{\bega{cl}\Hat W_i(s)\qquad &\hbox{if}\qquad t\leq \tau_\ve(\theta),\qquad 
\chi(\theta,\cdot)\Big|_{[0,t]}
~\simeq~\Hat \gamma_i\Big|_{[0,s]}\,,\\[4mm]
0\qquad &\hbox{if}\qquad t>\tau_\ve(\theta).\enda\right.\eeq
As proved in \cite{BSUN2}, the map $\ve\mapsto W^\ve(\theta,t)$ is nondecreasing for each $(\theta,t)$.
This leads to:

\begin{definition} Let $\chi:\Theta\times \R_+\mapsto\R^d$ be an irrigation plan 
	satisfying {\bf (A2)}.
	The {\bf weight function} $W=W(\theta,t)$ for $\chi$ is  defined as
	\bel{W}
	W(\theta,t)~\doteq~\sup_{\ve >0} ~W^\ve(\theta,t).\eeq
\end{definition}

Once we computed the weight functions on the irrigation plan $\chi$, its weighted irrigation cost $\E^{W,\alpha}$ is defined as follows:
\begin{definition} Let $f:\R_+\mapsto\R_+$ be a continuous function, satisfying
	all the assumptions in  {\bf (A1)}.
	Let $\chi $ be an irrigation plan satisfying {\bf (A2)} and let $W=W(\theta,t)$ be the 
	corresponding weight function, as in (\ref{W}). The 
	weighted cost $\E^{W,\alpha}$ for some $\alpha\in [0,1]$ is 
	\bel{10220}  \E^{W,\alpha}(\chi)~\doteq~\int_0^M \int_0^{\tau(\theta)} { (W(\theta,t))^\alpha \over m(\theta,t) } |\dot{\chi}(\theta, t)|\, dt\, d\theta\, . \eeq
\end{definition} 

In the special case where $\chi$ consists of only finitely many branches, let $W_i$ be the corresponding weight functions on the branch $\gamma_i:[a_i, b_i]\mapsto \R_+$, by applying the change of variable formula, we have the following identity for the weighted irrigation costs\cite{BSUN2}:
\bel{1a}
 \E^{W,\alpha}(\chi)~=~\sum_{i=1}^N \int_{a_i}^{b_i} [W_i(s)]^\alpha\, ds\, ,
\eeq 
where $N$ is the total number of branches.

\subsection{Lower semicontinuity of weighted cost}
In this section we recall the main  theorems on the lower semicontinuity of weighted irrigation cost,
proved  in \cite{BSUN2}.
Given  a sequence of irrigation plans $\chi_n: \Theta\times \R_+\mapsto \R^d$, we say that $\chi_n$ converges to $\chi$ pointwise if, for  every $\kappa>0$ and a.e.~$\theta\in \Theta$,
\bel{pconv}  \lim_{n\to \infty} \| \chi_n(\theta,\cdot)  - \chi(\theta,\cdot)  \|_{\mathbf{L}^\infty([0,\kappa])}~=~0. \eeq 

\begin{theorem}\label{t:lsc}
	
Let $(\chi_n)_{n\geq 1}$ be a sequence of irrigation plans, all satisfying {\bf (A2)}, pointwise converging
 to an irrigation plan $\chi$. Assume that the function $f$ satisfies {\bf (A1)}. Then 
\bel{lsc} \E^{W,\alpha}(\chi)~\leq~\liminf_{n\to \infty} \E^{W,\alpha}(\chi_n).   \eeq
	
\end{theorem}

Given a positive, bounded Radon measure $\mu$ on $\R^d$, we define the weighted irrigation cost $\I^{W,\alpha}(\mu)$ of $\mu$ as 
\bel{8}  \I^{W,\alpha}(\mu)~\doteq~\inf_\chi \E^{W,\alpha}(\chi), \eeq
where the infimum is taken over all irrigation plans for the measure $\mu$, and $\E^{W,\alpha}$ is defined as  in (\ref{10220}). 
By Theorem \ref{t:lsc}, if there is an irrigation plan for $\mu$ with finite weighted irrigation cost, then the infimum in (\ref{8}) is actually a minimum. That is, there exists an optimal irrigation plan $\chi^*$ of $\mu$, such that the weighted irrigation cost $\E^{W,\alpha}(\chi^*)$ is minimum among all admissible irrigation plans, and $\I^{W,\alpha}(\mu) = \E^{W,\alpha}(\chi^*)$.

The next result states the lower semicontinuity of the weighted irrigation cost, w.r.t.~weak convergence of measures.
For a proof,  see Theorem 6.2 in \cite{BSUN2}.
\begin{theorem}\label{t:meas}
	
	Let  $f$  satisfies {\bf (A1)}. Let $(\mu_n)_{n\geq 1}$ be a sequence of bounded positive Radon measures, with uniformly bounded supports, such  that weakly converges to some $\mu$. Then
	\bel{9} \I^{W,\alpha}(\mu)~\leq~\liminf_{n\to \infty} \I^{W,\alpha}(\mu_n).  \eeq
\end{theorem}

\section{Irrigability dimension}
\label{s:5}
\setcounter{equation}{0}

When $f=0$, $\alpha>1-{1\over d}$, 
it is well known
that  all measures with bounded support and finite mass in $\R^d$ are $\alpha$-irrigable \cite{BCM,DS}.
Here is a formal computation in this direction.
It is obtained by modifying the estimates at p.~113 of
\cite{BCM}.

Let $\mu$ be a probability measure that supported in $B(0,1)\subseteq \R^d$.
For each $j=1,\ldots,n$, let $\P_j$ be the set of centers of  balls of radius $r_j= 2^{-j}$ that cover $Supp(\mu)$.
In dimension $d$, we can assume that the  cardinality of this set is
$$\#\P_j~\leq~C \,2^{jd}$$
We can define a map $\gamma_j:\P_j\mapsto \P_{j-1}$ such that $$|x-\gamma_j(x)|~\leq ~3\cdot 2^{-j}$$
for every $x\in \P_j$, with $\P_0\doteq \{0 \}$.

Consider a probability measure $\mu_n$, supported on $\P_n$.
The cost of transporting this measure from $\P_n$ to another measure supported on $\P_{n-1}$ is 
\bel{ean}\bega{rl}
E^\alpha(\P_n, \P_{n-1})&\leq ~[\hbox{number of arcs}] \times [\hbox{flow}]^\alpha \times [\hbox{length}]\\[3mm]
&\leq
C\, 2^{nd} \cdot \left({1\over C\, 2^{nd}}\right)^\alpha\cdot 3\cdot 2^{-n}
~=~3 C^{1-\alpha} \cdot (2^{\alpha d-d+1})^{-n}.\enda\eeq
Notice that we are here considering the worst possible case, where we have the largest number of arcs and all arcs carry equal flow.

Summing over $j=1,2,\ldots,n$ we obtain that the total transportation cost is bounded by
\bel{Eser}E^\alpha~\leq~3 C^{1-\alpha} \cdot \sum_{j=1}^n (2^{\alpha d-d+1})^{-j}~\leq~
{3 C^{1-\alpha}\over  2^{\alpha d-d+1} -1}\,.\eeq
The series ~~$\sum_k 2^{[(d-1)-\alpha d](k+1)}$~~
converges provided that 
$$(d-1)-\alpha d~<~0,\qquad \hbox{hence}\qquad \alpha~>~1-{1\over d}$$

\begin{figure}[ht]
	\c{\includegraphics[width=14cm]{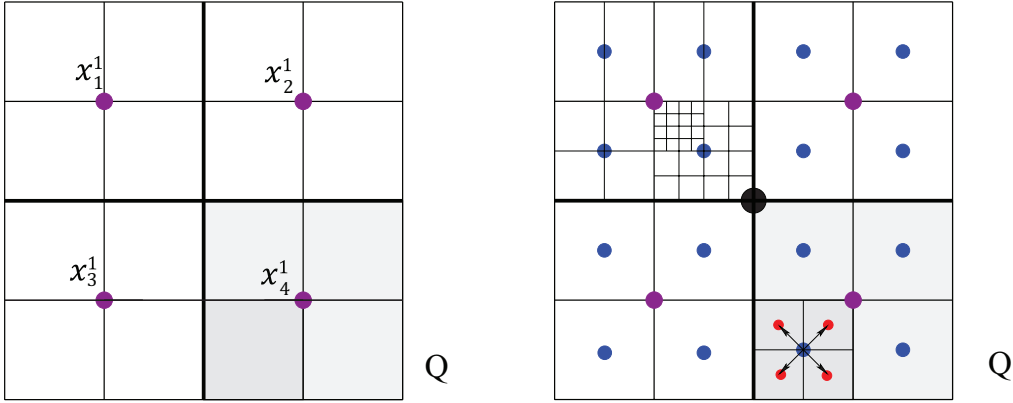}}
	\caption{\small Left: The dyadic approxmiated measure $\mu_1$ is supported on the four centers $x_1^1,\ldots,x^1_4$ of small cubes. Right: Dyadic approximated measures corresponding to a family of partitions into dyadic cubes in $\R^2$.     }
	\label{f:ir49}
\end{figure}

\v
To understand what happens in the case where weights are present, we first make an explicit computation
in the case of a dyadic irrigation plan \cite{BCM,X3}.
More precisely, as shown in Fig.~\ref{f:ir49}, we now assume

$\mu$ = Radon measure with total mass $M$,  concentrated on a cube $Q$ in $\R^d$. $Q$ is centered at the origin and with edge size $L>0$.

For each $n\geq 1$, we divide $Q$ into $2^{nd}$ smaller cubes of equal size, with edge size $L/2^n$. Take $\{Q_{i}^n \}_{i=1}^{2^{nd}}$ the set of all these closed smaller cubes, call $\P_n\doteq \{  x^n_i \}_{i=1}^{2^{nd}}$ the set of centers of these smaller cubes of edge size $L/2^n$.  For each $n\geq 1$, define the dyadic approximated measure $\mu_n$ 
\bel{1011} \mu_n~\doteq~\sum_{x^n_i\in \P_n}m_i^n\, \delta_{x_i^n},\eeq
where $\delta_{x_i^n}$ is the Dirac measure at $x_i^n$, and $m_i^n$ is determined as
$$\Hat Q_{i}^n~\doteq~Q_{i}^n\setminus \bigcup_{j<i} Q_{j}^n,\qquad m_i^n\doteq\mu(\Hat Q_{i}^n),\qquad \forall 1\leq i\leq 2^{nd}.$$

	It is not hard to show that $\mu_n$ weakly converges to $\mu$, see for example \cite{BCM,X3}. That is, for any bounded continuous function $\phi:\R^d \mapsto \R$, one has $\int \phi\, d\mu_n \to \int \phi \, d\mu. $ 
	For each $\mu_n$, we  construct  an irrigation plan $\chi_n$ as follows:
		\begi
	\item First, move the particles from the origin (center of $Q$) to the centers in $\P_1$, with $2^d$ straight paths connecting the origin and  the centers in $\P_1=\{x_1^1,x_2^1,\ldots,x^1_{2^d} \}$. Each path has length ${\sqrt{d}L\over 4}$, on the path that connecting $x^1_i, 1\leq i\leq 2^d$, the multiplicity is constant $m^1_i$.
	\item By induction, at the level $k, 1< k\leq n$, for the particles arriving at each center $x_i^{k-1}$ in $\P_{k-1}$, where $x_i^{k-1}$ is the center of the cube $Q^{k-1}_i$, we transport them to the $2^d$ neighboring centers in $\P_k$, which are all contained in the cube $Q_i^{k-1}$. Without loss of generality, fixed $x_i^{k-1}$ in $\P_{k-1}$, let $\{x^k_1,\ldots,x^k_{2^d} \}$ be the $2^d$ neighboring centers around $x_i^{k-1}$. For each $x_j^k,1\leq j\leq 2^d$, we build a straight path connecting $x_i^{k-1}$ to $x^k_j$, with length ${\sqrt{d}L\over 2^{k+1}}$ and constant multiplicity $m_j^k$.

	\endi
	\begin{figure}[hb]

	\c{\includegraphics[width=14cm]{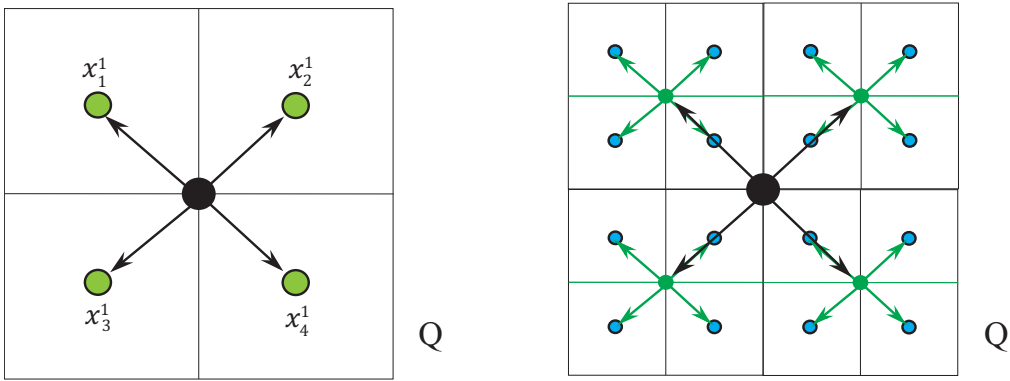}}
	\caption{\small The dyadic irrigation plans in $\R^2$. Left: The dyadic irrigation plan $\chi_1$. The multiplicity on each branch equals to the mass on the terminal point. Right: The dyadic irrigation plan $\chi_2$. The particles are first transported to the 4 centers in $\P_1$, then on each center in $\P_1$, the particles are transported to the neighboring 4 centers in $\P_2$. }
	\label{f:dir7}
\end{figure}
Since the dyadic measure $\mu_n$ is supported on the centers in $\P_n$, after $n$ steps we build an irrigation plan for $\mu_n$, which we call the {\bf dyadic irrigation plan} $\chi_n$.

For example, in the case $\R^2$, Fig. \ref{f:dir7} shows two dyadic irrigation plans constructed by the preceding procedure.

Given an irrigation plan with finite branches as in Section \ref{finite}, consider the case $f(z) \doteq cz^\beta$, with some constant $c>0, 0<\beta<1$. It is readily to check that $f$ satisfies {\bf (A1)}.  With the notions in Section \ref{finite}, consider a measure $\mu$ consisting of finitely many point masses $m_i\geq 0$ located at points $P_i$, where $P_i$ is the ending node of branch $\gamma_i(s):[0,\ell_i]\mapsto \R^d$. In this case, the multiplicity function on each branch is constant. Then the computation of weights (\ref{qs1})-(\ref{qs3}) becomes 
\bel{qs9}
\bega{c}\ds
W_i(s)~=~\left(\ov W_i^{1-\beta} + c(1-\beta)(\ell_i -s)\right)^{{1\over 1-\beta}}\, ,\\[3mm]
\ov W_i ~=~m_i + \sum_{j\in \O(i)} \left(\ov W_j^{1-\beta} + c(1-\beta)\ell_j \right)^{{1\over 1-\beta}}\, .
\enda\eeq
If $\O(i)=\emptyset$, that is $i\in \I_1$, from (\ref{qs9}) we have $\ov W_i = m_i$.

	 The following two lemmas proved that under suitable conditions, the weighted irrigation costs of the dyadic irrigation plans $\{\chi_n \}_{n\geq 1}$ are uniformly bounded. Utilizing this fact and Theorem \ref{t:meas}, since the dyadic approximated measures $\mu_n$ weakly converges to $\mu$, we can conclude the irrigability of $\mu$ with weighted cost.

	To fix the ideas, we first consider the case that $\mu$ is the Lebesgue measure on the unit cube $Q$.
		\begin{lemma}\label{l:41}	
		Suppose $1>\beta > 1 - {1\over d}$, $1\geq \alpha > 1- {1\over d}$, while 
		$\mu$ is the Lebesgue measure on the unit cube $Q$ in $\R^d$.  Then, in the dyadic irrigation plans $\chi_n$, the weight function $W^n$ remains uniformly bounded on all branches.	
		Moreover, the irrigation cost $\E^{W,\alpha}(\chi_n)$ is uniformly bounded. That is, there exists an uniform constant $C>0$, such that for all dyadic irrigation plan $\chi_n$,
		\bel{uniw} W^n~\leq~C,\qquad \E^{W,\alpha}(\chi_n)~\leq~C \, .\eeq
	\end{lemma}

	{\bf Proof. } For the dyadic irrigation plan $\chi_n$, since each dyadic irrigation  plan has finite branches and $\mu_n$ is supported on the centers in $\P_n$, we can use formula (\ref{qs9}) to compute the weights $W^n$. We start from the centers in $\P_n$.
	
	{\bf 1.} From $\P_n$ to $\P_{n-1}$,  for each $x_i^{n-1}\in \P_{n-1}$, by the construction of the dyadic irrigation plan $\chi_n$, there are straight paths connecting $x_i^{n-1}$ to the $2^d$ neighboring centers in $\P_n$. 
	 Since $\mu$ is the Lebesgue measure on unit cube $Q$, 
	mass on each center in $\P_n$ is ${1\over 2^{nd}}$. The branches connecting $x^{n-1}_i$ to centers in $\P_n$ are identical, with branch length $\sqrt{d}/2^{n+1}$ and constant multiplicity $1/2^{nd}$. We only need  to  compute the weight on one such branch, and write it as $W_n^n$, where the superindex $n$ means it is the weight for irrigation plan $\chi_n$, and the subindex $n$ means from $\P_n$ to $\P_{n-1}$. 
	
	By formula (\ref{qs9}), for $s\in\, [0,\, {\sqrt{d}\over 2^{n+1}}]$,
	\bel{w1}W^n_n(s)~=~\Big(({1\over 2^{nd}})^{1 -\beta} + c(1 - \beta)({\sqrt{d}\over 2^{n+ 1}} - s)\Big)^{1\over 1 -\beta}\, ,   \eeq
	\bel{w2}W^n_n(0)~=~\Big(({1\over 2^{nd}})^{1 -\beta} + c(1 - \beta){\sqrt{d}\over 2^{n+ 1}}\Big)^{1\over 1 -\beta}.\eeq
	\v
	{\bf 2.} From $\P_{n-1}$ to $\P_{n-2}$, using formula (\ref{qs9}), on each branch we need to first compute the weights  $\overline{W}_{n-1}^n$ at the tip. For the dyadic approximated measure $\mu_n$, it is supported on $\P_n$, thus the mass on each center in $\P_k, k\neq n$ is $0$. Since  each center in $\P_{n-2}$ connects $2^d$ identical centers in $\P_{n-1}$, we therefore have
	\bel{w3}\overline{W}^n_{n-1}~=~2^dW_n^n(0)~=~2^d\Big(({1\over 2^{nd}})^{1 -\beta} + c(1 - \beta){\sqrt{d}\over 2^{n+ 1}}\Big)^{1\over 1 -\beta}.\eeq
	Each branch between $\P_{n-2}$ and $\P_{n-1}$ has length ${\sqrt{d}\over 2^n}$. By formula (\ref{qs9}), for $ s \in \, [0,\, {\sqrt{d}\over 2^n} ]$, 
	\bel{w4}W^n_{n-1}(s)~=~\Big(2^{d(1 -\beta)}\Big[({1\over 2^{nd}})^{1 -\beta} + c(1 - \beta){\sqrt{d}\over 2^{n +1}}\Big] + c(1 -\beta)({\sqrt{d}\over 2^n} - s)\Big)^{1\over 1 -\beta}\, , \eeq
	\bel{w5}W^n_{n- 1}(0)~=~\Big(({1\over 2^{(n-1) d}})^{1 -\beta} + c(1 -\beta) \Big[{\sqrt{d}\over 2^{n + 1 - d(1 -\beta)}} + {\sqrt{d}\over 2^n} \Big] \Big)^{1\over 1-\beta}\, .\eeq
	\v
	{\bf 3.} From $\P_{n-2}$ to $\P_{n-3}$. At the tip of each branch,
	\bel{w6} \overline{W}_{n-2}^n~=~2^d W^n_{n - 1}(0). \eeq
	The branches have length ${\sqrt{d}\over 2^{n-1}}$, by formula (\ref{qs9}), for $ s\in \,  [ 0,\, {\sqrt{d}\over 2^{n -1}}]$,
	\bel{w7}\bega{rl}\ds W^n_{n - 2}(s)~=&\ds \Big( 2^{d(1 -\beta)}\left[({1\over 2^{(n-1) d}})^{1 -\beta} +  c(1 -\beta) \Big[{\sqrt{d}\over 2^{n + 1 - d(1 -\beta)}} + {\sqrt{d}\over 2^n} \Big]  \right]\\[3mm]
	&\ds\qquad+~ c(1 -\beta)({\sqrt{d}\over 2^{n - 1}}-s) \Big)^{1\over 1 - \beta}\\[3mm]
	~=&\ds \Big(({1\over 2^{(n-2) d}})^{1 -\beta} + c(1 -\beta)  \Big[{\sqrt{d}\over 2^{n + 1 - 2d(1 - \beta)}} + {\sqrt{d}\over 2^{n - d(1 -\beta)}}+~  ({\sqrt{d} \over 2^{n -1}}-s)\Big] \Big)^{1\over 1 -\beta}\, ,
	\enda\eeq
	$$W_{n-2}^n(0)~=~\Big(({1\over 2^{(n-2) d}})^{1 -\beta} +c(1 -\beta)  \Big[{\sqrt{d}\over 2^{n + 1 - 2d(1 - \beta)}} + {\sqrt{d}\over 2^{n - d(1 -\beta)}}+~  {\sqrt{d} \over 2^{n -1}}\Big] \Big)^{1\over 1 -\beta}.$$
	\v
	{\bf 4.} From $\P_{n - k}$ to $\P_{n - k -1}$, each branch has length ${\sqrt{d}\over 2^{n+1-k}}$. Similarly for $s\in \, [0,\, {\sqrt{d}\over 2^{n+1-k}}]$,
	$$W^n_{n - k}(s)=\Big(({1\over 2^{(n-k) d}})^{1 -\beta} + {c(1-\beta)\sqrt{d}\over 2^{n+1}}\sum_{j=1}^{k } 2^{(k-j)+ jd(1 - \beta)} + c(1 -\beta)({\sqrt{d}\over 2^{n+1-k}} - s)\Big)^{1\over 1 -\beta},$$
	\bel{w9}\bega{rl}
	\ds W_{n-k}^n(0)~=&\ds\Big(({1\over 2^{(n-k) d}})^{1 -\beta} + {c(1-\beta)\sqrt{d}\over 2^{n+1}}\sum_{j=0}^{k } 2^{(k-j)+ jd(1 - \beta)}\Big)^{1\over 1 -\beta}\\[3mm]
	~=&\ds\Big(({1\over 2^{(n-k) d}})^{1 -\beta} + {c(1 -\beta)\sqrt{d} \over 2^{n+1 - k}}\sum_{j=0}^k {1\over 2^{[1 -d(1 - \beta)]j}}\Big)^{1\over 1 -\beta}.
	\enda\eeq
	\v
	{\bf 5.} Since $W^n_{n-k}(s)\leq W^n_{n-k}(0)$, we only need to estimate $W^n_{n - k}(0)$, for each $0\leq k\leq n-1.$ When $\beta > 1 - 1/d$, one has $1 - d(1 -\beta) >0$. By (\ref{w9}), for each $k$,
	\bel{w11} W^n_{n - k}(0)~\leq~\Big(({1\over 2^{(n-k) d}})^{1 -\beta} + {c(1 -\beta)\sqrt{d}\over 2^{n + 1 - k}} \cdot {1\over 1 - ({1\over 2})^{1 - d(1-\beta)} }\Big)^{1\over 1 -\beta}.\eeq
	Therefore, we have an uniform bound for the weight function
	\bel{w18} W^n
	~\leq~\left(1 + {c(1-\beta)\sqrt{d}\over1- ({1\over 2})^{1 - d(1-\beta)} } \right)^{1\over 1-\beta},\eeq
	which is independent of $n$.
	\v
	{\bf 6.} We now estimate the irrigation cost $\E^{W,\alpha}(\chi_n)$ by the formula (\ref{1a}). Fixed the dyadic irrigation plan $\chi_n$, call $E_n^n$ the  cost from $\P_n$ to $\P_{n-1}$.
	There are $2^{nd}$ branches from centers in $\P_n$ to centers in $\P_{n-1}$. On each branch, the weight $W^n_n$ is given by (\ref{w1}). Therefore,
	\bel{c1}\bega{rl}\ds E^n_n&\ds=~2^{nd}\int_0^{\sqrt{d}\over 2^{n+1}}\Big(({1\over 2^{nd}})^{1 -\beta}+ c(1 -\beta)({\sqrt{d}\over 2^{n+1}} - s)\Big)^{\alpha\over 1 -\beta}\, ds\\[3mm]
	&\ds=~{2^{nd}\over c(1 + \alpha - \beta)}\left(\left[({1\over 2^{nd}})^{1 -\beta} + c(1 -\beta){\sqrt{d}\over 2^{n+1}}\right]^{1+\alpha-\beta\over 1-\beta} - \left[({1\over 2^{nd}})^{1-\beta}\right]^{1 + \alpha -\beta\over 1-\beta}\right).
	\enda\eeq
	Similarly, denote $E^n_{n-k}$ the  cost from $\P_{n-k}$ to $\P_{n-k-1}$. There are  $2^{(n-k)d}$ branches from centers in $\P_{n-k}$ to centers in $\P_{n-k-1}$. 
	\bel{c2}\bega{rl}\ds E^n_{n -k}
	&\ds =~2^{(n - k)d}\int_0^{\sqrt{d}\over 2^{n+1 - k}}\Big(\left(\overline{W}^n_{n-k}\right)^{1-\beta} + c(1 -\beta)({\sqrt{d}\over 2^{n+ 1 - k}} - s)\Big)^{\alpha\over 1-\beta}\, ds\\[3mm]
	&\ds =~{2^{(n - k)d}\over c( 1 + \alpha - \beta)}\Big(\Big[\left(\overline{W}^n_{n-k}\right)^{1-\beta} + c(1 -\beta){\sqrt{d}\over 2^{n +1 - k }}\Big]^{1 + \alpha -\beta\over 1 -\beta} - \left(\overline{W}^n_{n-k}\right)^{1+\alpha -\beta}\Big).
	\enda
	\eeq
	In the following, we use the same $C$ to denote different constants which only depend on $c,\alpha,\beta$ and the dimension $d$. From (\ref{w11}) and the fact that $(1 - \beta)d <1$, for each $n$ and $k$,
	\bel{rati}\left(\overline{W}^n_{n - k}\right)^{1 -\beta} + c(1 -\beta){\sqrt{d}\over 2^{n +1 - k }}~\leq~{C\over 2^{(n-k)(1 -\beta)d}}\, .\eeq
	Consider $x, y\geq 0$, 
	\bel{c3}g(x , y)~\doteq~(x + y)^{1 +\alpha -\beta\over 1 -\beta} - x^{1 + \alpha -\beta\over 1-\beta}\, ,\qquad \,  x+ y \leq {C\over 2^{(n - k)(1 -\beta)d}}\,.\eeq
	then, by a first order Taylor expansion,
	\bel{c4}g(x , y)~\leq~{1 +\alpha-\beta\over 1-\beta}\left({C\over 2^{(n-k)(1-\beta)d}} \right)^{{\alpha\over 1-\beta }}\cdot y~\leq~{C y\over 2^{(n-k)\alpha d}}.\eeq
	Applying (\ref{rati}) and (\ref{c4}) in (\ref{c2}), we obtain
	\bel{c5}E^n_{n - k}~\leq~2^{(n-k)d}{C \over 2^{(n-k)(\alpha d + 1)}	}~=~{C\over 2^{(n-k)[(\alpha - 1)d + 1]}}\, .\eeq
	When $\alpha > 1 -1/d$, one has $(\alpha - 1)d + 1 > 0$. Then by (\ref{c5}),
		\bel{c6}\E^{W,\alpha}(\chi_n)~=~\sum_{k =0}^{n-1}E^n_{n -k}~\leq~\sum_{k=0}^{n-1}{C\over 2^{(n-k)[(\alpha-1)d + 1]} }~\leq~{C\over 1 -  ({1\over 2})^{(\alpha - 1)d +1}},\eeq
	where $C$ is some constant independent of $n$.
	Combining the estimates (\ref{w18}) and (\ref{c6}), 
	we obtain the existence of a constant $C$, independent of $n$, such that (\ref{uniw}) holds.
		\endproof

	Under the same conditons on $\alpha, \beta$, this uniform boundedness result holds for general positive, finite Radon measures.
	\begin{lemma}\label{l:42}
		Suppose $ 1 - {1\over d} < \beta < 1$, $ 1 -{1\over d} < \alpha \leq 1$. $\mu$ is a finite measure on the cube $Q$ with edge size $L$ in $\R^d$, denote $M$ the total mass of $\mu$. Then in the dyadic irrigation plan $\chi_n$, the weight function $W^n$ on each branch remains uniformly bounded, 
		\bel{uniw1} W^n~\leq~C\Big(M^{1 -\beta} + L\Big)^{1\over 1-\beta}\, ,\eeq
		Moreover, the irrigation cost $\E^{W,\alpha}(\chi_n)$ is uniformly bounded, namely
		\bel{univ1}\E^{W,\alpha}(\chi_n)~\leq~C\Big(M^\alpha L + L^{1 + {\alpha\over 1 - \beta}}\Big) \eeq
		where $C$ is some constant independent of $n$.
	\end{lemma}
	{\bf Proof.}  For the dyadic irrigation plan $\chi_n$, to compute the weights $W^n$, we start from the centers in $\P_n$.

	{\bf 1.} From $\P_n$ to $\P_{n-1}$. Let $m_i^n$ be the mass of $\mu_n$ on the center $x_i^n$ in $\P_n$. On the branch from $x^n_i$ to any center in $\P_{n-1}$, the arc-length of the branch is ${\sqrt{d}L\over 2^{n+1}}$ and the multiplicity is constant $m_i^n$. Let $W^n_{n,i}$ be the corresponding weights, where the superindex $n$ means we consider the weight function on irrigation plan $\chi_n$, the subindex $(n,i)$ means we consider the weight on the $i$-th branch from $\P_n$ to $\P_{n-1}$.  Then  by formula (\ref{qs9}), for $s\in\, [ 0,\,   {\sqrt{d}L\over 2^{n + 1}}  ]$,
	\bel{nw1}W_{n,i}^n(s)~=~\Big(  \left(m_i^n\right)^{1 -\beta} + c(1 -\beta)({\sqrt{d}L\over 2^{n+1}} - s) \Big)^{1\over 1 -\beta}\,, \eeq
	\bel{nw2}W_{n,i}^n(0)~=~\Big(\left(m_i^n\right)^{1 -\beta} +c(1 -\beta){\sqrt{d}L\over 2^{n+1}} \Big)^{1\over 1 -\beta}\, .\eeq
	\v
	{\bf 2.} From $\P_{n-1}$ to $\P_{n-2}$. For each center $x^{n-1}_i$ in $\P_{n-1}$, to compute the weight $W^n_{n-1,i}$ from $x^{n-1}_i$ to any center in $\P_{n-2}$, we first estimate  $\overline{W}^n_{n-1,i}$. Each  $x^{n-1}_i$ in $\P_{n-1}$ connects $2^d$ nearby centers in $\P_n$. By (\ref{qs9}) and (\ref{nw2}) one has,
	\bel{nw3}\overline{W}^n_{n-1,i}~=~\sum_{j\in \O(i)} W^n_{n, j}(0)~=~\sum_{j\in\O(i)}\Big(\left(m_j^n\right)^{1 -\beta} + c(1 -\beta) {\sqrt{d}L\over 2^{n+1}}\Big)^{1\over 1 -\beta}\,.\eeq
	Notice for fixed $b\geq 0$, $g(x)\doteq (x^{1 -\beta} + b)^{1\over 1 -\beta}$ is a concave function of $x$ on $\R_+$. Thus for any $N$,
	\bel{cova}{1\over N}\sum_{j=1}^N (x_j^{1 -\beta} + b)^{1\over 1-\beta}~\leq~\Big(({\sum_{j=1}^N x_j\over N})^{1-\beta} + b\Big)^{1\over 1-\beta}\, .\eeq
	For each $i$, the cardinality of the set $\O(i)$ in (\ref{nw3}) is $2^d$. From (\ref{nw3})-(\ref{cova}),
	\bel{nw4}\overline{W}^n_{n-1,i}~\leq~2^d\Big[\Big({\sum_{j\in\O(i)} m^n_j\over 2^d}\Big)^{1 -\beta} + c(1 -\beta){\sqrt{d}L\over 2^{n+1}}\Big]^{1\over 1 -\beta}\,.\eeq
	Each branch from $x_i^{n-1}$ to $\P_{n-2}$ has length ${\sqrt{d}L\over 2^n}$. By the formula (\ref{qs9}), for $s\in \, [0,\, {\sqrt{d}L\over 2^n}]$,
	\bel{nw5}\bega{c}\ds W^n_{n-1, i}(s)~ =~\Big(\left(\overline{W}_{n-1, i}^n\right)^{1 -\beta} + c(1 -\beta)({\sqrt{d}L\over 2^{n}} - s) \Big)^{1\over 1-\beta}\\[3mm]
	\ds \leq~\Big( (\sum_{j\in \O(i)} m_j^n)^{1 -\beta} + c(1-\beta)\Big[ {\sqrt{d}L\over2^{n+ 1 -d(1-\beta)}} + ({\sqrt{d}L\over2^n } - s) \Big] \Big)^{1\over 1-\beta}\, ,
	\enda\eeq
	\bel{nw6}W^n_{n-1, i}(0)~\leq~\Big( (\sum_{j\in \O(i)} m^n_j)^{1 -\beta} + c(1-\beta) \Big[ {\sqrt{d}L\over2^{n+ 1 -d(1-\beta)}} + {\sqrt{d}L\over2^n } \Big] \Big)^{1\over 1-\beta}\, .\eeq
	\v
	{\bf 3.} From $\P_{n-2}$ to $\P_{n - 3}$. For each center $x^{n-2}_i$ in $\P_{n-2}$,
	according to (\ref{nw6}),
	\bel{nw7}\bega{l}\ds \overline{W}^n_{n - 2, i}~=~\sum_{k\in \O(i)} W^n_{n - 1, k}(0)\\[3mm]
	\ds \leq~\sum_{k\in \O(i)}\Big( (\sum_{j\in\O(k)} m_j^n)^{1 -\beta} + c(1-\beta) \Big[ {\sqrt{d}L\over2^{n+ 1 -d(1-\beta)}} + {\sqrt{d}L\over2^n } \Big]  \Big)^{1\over 1-\beta}\,.\enda\eeq
	Using the concavity inequality (\ref{cova}), 
	\bel{nw8}  \overline{W}^n_{n - 2, i}~\leq~2^d\left[\Big( {\sum_{k\in \O(i),j\in\O(k)}m_j^n\over 2^d}\Big)^{1 -\beta} + c(1-\beta) \Big[ {\sqrt{d}L\over2^{n+ 1 -d(1-\beta)}} + {\sqrt{d}L\over2^n } \Big]  \right]^{1\over 1-\beta}\, .\eeq
	In the following, for each center $x^k_i$ in $\P_k$, if there is a concatenated path from  $x^k_i$ to center $x^n_j\in \P_n$ in the dyadic irrigation plan $\chi_n$, we  say $i\prec j$.  With this notation, (\ref{nw8}) can be written as 
	\bel{1012}  \overline{W}^n_{n - 2, i}~\leq~2^d\left[\Big( {\sum_{i\prec j}m_j^n\over 2^d}\Big)^{1 -\beta} + c(1-\beta) \Big[ {\sqrt{d}L\over2^{n+ 1 -d(1-\beta)}} + {\sqrt{d}L\over2^n } \Big]   \right]^{1\over 1-\beta}\, .  \eeq
	Each branch from $x_i^{n-2}$ to $\P_{n-3}$ has length ${\sqrt{d}L\over 2^{n-1}}$. By the formula (\ref{qs9}), for $s\in \, [0,\, {\sqrt{d}L\over 2^{n - 1}}]$,
	\bel{nw9}\bega{l}\ds W^n_{n-2, i}(s)~=~ \Big(\left(\overline{W}^n_{n -2, i}\right)^{1 -\beta} + c(1-\beta)({\sqrt{d}L\over 2^{n-1}} - s) \Big)^{1\over 1-\beta}\qquad\qquad\\[3mm]
	\ds\leq~\Big( (\sum_{i\prec j}m^n_j )^{1 -\beta} + c(1-\beta)\Big[{\sqrt{d}L\over 2^{n+ 1 - 2d(1 -\beta)}} +{\sqrt{d}L\over 2^{n - d(1 -\beta)}} + ({\sqrt{d}L\over 2^{n- 1}} -s)\Big]\Big)^{1\over 1-\beta}.
	\enda\eeq
	\bel{nw10}W^n_{n-2, i}(0)~\leq~\Big((\sum_{i\prec j}m^n_j)^{1 -\beta} + c(1-\beta)\Big[{\sqrt{d}L\over 2^{n+ 1 - 2d(1 -\beta)}} +{\sqrt{d}L\over 2^{n - d(1 -\beta)}} + {\sqrt{d}L\over 2^{n- 1}}\Big]\Big)^{1\over 1-\beta} \, .\eeq
	\v
	
	{\bf 4.} From $\P_{n - k}$ to $\P_{n - k -1}$ . Similarly we have,
	\bel{nwa}\overline{W}^n_{n-k,i}~\leq~2^d\Big(({\sum_{i\prec j}m^n_j\over 2^d})^{1-\beta} + {c(1-\beta)\sqrt{d}L\over 2^{n+ 2 -k} }\sum_{l=0}^{k-1} {1\over 2^{[1- d(1 -\beta)]l}} \Big)^{1\over 1-\beta}\, ,\eeq
	\bel{nw11} W^n_{n - k, i}(s)\leq\Big((\sum_{i\prec j}m^n_j)^{1 -\beta} +{c(1-\beta)\sqrt{d}L\over 2^{n+ 1 -k} }\sum_{l=1}^{k} {1\over 2^{[1- d(1 -\beta)]l}} + c(1-\beta)({\sqrt{d}L\over2^{n +1 - k}}-s)\Big)^{1\over 1-\beta}\eeq
	\bel{nw12}
	W^n_{n-k,i}(0)~\leq~\Big((\sum_{i\prec j}m^n_j)^{1 -\beta} + {c(1 -\beta)\sqrt{d}L\over 2^{n+1 - k}}\sum_{l =0}^k {1\over 2^{[1 -d(1 -\beta)]l}}\Big)^{1\over 1-\beta}\, .
	\eeq
	\v
	
	{\bf 5.} Since $W^n_{n-k,i}(s)\leq W^n_{n-k,i}(0)$, we only need to estimate $W^n_{n-k,i}(0)$, for each $0\leq k\leq n-1, 1\leq i\leq 2^{d(n - k)}$. When $\beta > 1 -1/d$, one has $1 - d(1 -\beta)>0$. From formula (\ref{nw12}), 
	\bel{nw13}W^n_{n-k,i}(0)~\leq~\Big((\sum_{i\prec j}m^n_j)^{1 -\beta} + {c(1-\beta)\sqrt{d}L\over 2^{n+1 - k}}{1 \over 1 -{1\over 2^{1- d(1 -\beta)} }} \Big)^{1\over 1-\beta}.\eeq
	Since $\sum_{i\prec j}m^n_j\leq M $, if denote $W^n$ the weights on dyadic irrigation plan $\chi_n$, from (\ref{nw13}) 
	there is an uniform bound for the weight function 
	\bel{nnw}W^n~\leq~\Big(M^{1 -\beta} + {c(1 -\beta)\sqrt{d}L\over 1 - {1\over 2^{1- d(1 -\beta)} }}\Big)^{1\over 1-\beta}~\leq~C\Big(M^{1 -\beta} + L \Big)^{1\over 1-\beta}\eeq
	where we use the same $C$ to denote all constants independent of $n$. This completes the proof of (\ref{uniw1}).
	\v
	
	{\bf 6.} We now estimate the irrigation cost $\E^{W,\alpha}(\chi_n)$ by formula (\ref{1a}). In the dyadic irrigation plan $\chi_n$, let $E^n_n$ be the cost from $\P_n$ to $\P_{n- 1}$, by (\ref{nw1}),
	\bel{e1}\bega{rl}\ds E^n_n&\ds=~\sum_{x^n_i\in \P_n} \int_0^{{\sqrt{d}L\over 2^{n+1}}} \left(W^n_{n,i}(s)\right)^\alpha\, ds\\[3mm]
	&\ds=~\sum_{x^n_i\in \P_n} \int_0^{{\sqrt{d}L\over 2^{n+1}}} \Big(\left(m^n_i\right)^{1 -\beta} + c(1 -\beta)({\sqrt{d}L\over 2^{n + 1}} - s) \Big)^{\alpha\over 1 -\beta}\, ds\,.\\[3mm]
\enda\eeq
	Similarly, denote $E^n_{n-k}$ the  cost from $\P_{n- k}$ to $\P_{n- k - 1}$, 
	\bel{e2}\bega{l}\ds E^n_{n - k}~=~\sum_{x^{n-k}_i\in \P_{n - k}}\int_0^{{\sqrt{d}L\over 2^{n+ 1 -k}}} \left(W_{n-k,i}^n (s)\right)^\alpha\, ds\enda\eeq
	From (\ref{nw12}) and the non-decreasing of $W^n_{n-k,i}(s)$,
	\bel{e3}\bega{rl}\ds E^n_{n- k}&\ds\leq~\sum_{x_i^{n-k} \in \P_{n - k}}{\sqrt{d}L\over 2^{n+1 - k}}\Big((\sum_{i\prec j}m^n_j)^{1 -\beta} + {c(1-\beta)\sqrt{d}L\over 2^{n+1 - k} (1 -{1\over 2^{1- d(1 -\beta)}})} \Big)^{\alpha\over 1-\beta}\\[3mm]
	&\ds\leq~\sum_{x_i^{n-k}\in \P_{n - k}}\Big[{CL(\sum_{i\prec j} m^n_j)^{\alpha}\over 2^{n+ 1 -k}} + {C L^{1 + {\alpha\over 1-\beta}}\over 2^{(n+1 - k)({\alpha\over 1-\beta}+ 1 )} }\Big]~\doteq~I_{n-k} + J_{n-k}
	\enda\eeq
	where $C$ is some constant that only depends on $\alpha,\beta, c$ and on the dimension $d$.
	The cardinality of  $\P_{n - k}$ is $2^{(n - k)d}$. Therefore
	\bel{e4}J_{n-k}~\doteq~\sum_{x_i^{n-k}\in \P_{n - k}} {C L^{1 + {\alpha\over 1-\beta}}\over 2^{(n+1 - k)({\alpha\over 1-\beta}+ 1 )} }~\leq~  {CL^{1 + {\alpha\over 1-\beta}}\over 2^{(n-k)(1 + {\alpha\over 1-\beta} -d)}}\, .\eeq
	On the other hand, $1\geq \alpha>0$, by elementary concavity inequality,
	\bel{e5}\bega{rl}\ds I_{n-k}&\ds\doteq~\sum_{x_i^{n-k}\in \P_{n - k}} {CL(\sum_{i\prec j} m^n_j)^{\alpha}\over 2^{n+ 1 -k}}  ~\leq~2^{(n-k)d} \Big({\sum_{x_i^{n-k}\in \P_{n-k} }\sum_{i\prec j} m_j^n  \over 2^{(n-k)d} } \Big)^\alpha {CL\over 2^{n+1 - k}}\\[3mm]
	&\ds\leq~{CM^\alpha L \over  2^{(n-k) [ 1 - d(1-\alpha)  ]}  }.
	\enda\eeq 
	When $1\geq \alpha> 1-1/d$ and $1>\beta> 1-1/d$, one has
	\bel{e6}1 - d(1 -\alpha)~>~0\, ,\qquad 1+{\alpha\over 1-\beta} -d~>~0\, .\eeq
	Therefore, using (\ref{e3})-(\ref{e5}),
	\bel{e7}\bega{rl}\ds\E^{W,\alpha}(\chi_n)&\ds=~\sum_{k=0}^{n-1} E^n_{n - k}~\leq~\sum_{k=0}^{n-1} \left[ I_{n-k} + J_{n-k} \right]\\[3mm]
	&\ds\leq~C\sum_{j=0}^n\Big[{LM^\alpha\over 2^{[1 - d(1-\alpha)]j}} + {L^{1 + {\alpha\over 1 -\beta}}\over 2^{(1 +\frac{\alpha}{1-\beta} - d) j}}\Big]\\[3mm]
	&\ds\leq~C\Big(LM^\alpha + L^{ 1 + {\alpha\over 1 - \beta}}\Big)
	\enda\eeq
	where $C$ is some constant independent of $n$. This completes the proof of (\ref{univ1}).
		\endproof
\v
By the previous results, when $f(z)\doteq c z^\beta$ in (\ref{qs1})-(\ref{qs3}), with the conditions in Lemma \ref{l:42}, we have the uniform bounds (\ref{univ1}) for the dyadic irrigation plan sequence $\{\chi_n\}_{n\geq 1}$. Since each $\chi_n$ is an admissible irrigation plan for $\mu_n$, by the definition (\ref{8}), we have a uniform bound on all
the irrigation costs $\I^{W,\alpha}(\mu_n)$, $ n\geq 1$. By the weak convergence  $\mu_n\wto\mu$ and
the lower semicontinuity  of the irrigation cost, stated in Theorem \ref{t:meas}, 
we conclude $\I^{W,\alpha}(\mu)<+\infty$. 
 
By a comparison argument we can now prove the irrigability
for a wide class of functions $f$ and measures $\mu$, with the weighted irrigation cost $\I^{W,\alpha}$ in (\ref{8}).

\begin{theorem}\label{t:i}
	Let $\mu$ be a positive, bounded Radon measure  in $\R^d$, with total mass $M>0$ and  supported in the cube $Q$ of edge size $L>0$.  Assume $\alpha > 1 - {1\over d}$, $f$ satisfies {\bf (A1)} and
	\bel{fp4}
	\limsup_{z\to 0^+} ~z^{-\beta} f(z)~<~+\infty\eeq
	for some $1>\beta>1-{1\over d}$.   Then $\I^{W,\alpha}(\mu)<+\infty$.
\end{theorem}

{\bf Proof.} 
	The assumptions (\ref{fp4}) and {\bf (A1)} together imply that
	\bel{fp5}\bega{cl}\ds f(z)~\leq~c z^\beta&\ds\qquad\qquad\forall z\in [0, 2^{1\over 1-\beta}z_0],\\[3mm] \ds f(z)~\leq~cz&\ds\qquad\qquad\forall z\in [z_0,\infty),\enda\eeq
with some constants $c,z_0 >0$. We will prove that the weighted irrigation costs of the dyadic approximated measures $\mu_n$, defined as in (\ref{1011}), are uniformly bounded.  Since $\mu_n$ weakly converges to $\mu$, by Theorem \ref{t:meas}, this uniform bound implies the boundedness of $\I^{W,\alpha}(\mu).$

 It suffices to prove the uniform bound for  dyadic approximated measures $\mu_n = \sum_{x_i^n\in \P_n} m_i^n \delta_{x_i^n}$ with $n\geq n_0$, where $n_0$ is some fixed integer. Choose $n_0$ large enough such that in (\ref{nw13}),
\bel{1020}  {c(1-\beta)\sqrt{d}L\over 2^{n_0 }\cdot (1 - {1\over 2^{1 - d(1-\beta)}})}~<~z_0^{1-\beta}. \eeq In the following, we construct the  irrigation plan for $\mu_n$ with uniformly bounded weighted cost.

{\bf 1.} Consider first from $\P_n$ to $\P_{n-1}$. For those $x_i^n$ such that $m_i^n\geq z_0$, we transport the particles at $x_i^n$ along a straight path directly to the origin. Let $\S_n$ be the set of all such paths. For each path in $\S_n$, the multiplicity is larger than $z_0$ and bounded by $M$. The length of path is
 bounded by $\sqrt{d}L$. Let $W(t)$ be the  weight function on these paths, then clearly $W(t)\geq z_0$. By formula (\ref{qs1})-(\ref{qs3}) and (\ref{fp5}) the weight satisfies
\bel{1021} W(t)~\leq~\int_t^{\sqrt{d}L} f(W(s))\, ds + M ~\leq~\int_t^{\sqrt{d}L} c W(s)\, ds + M~\leq~e^{c\sqrt{d}L}M\, .\eeq  
On the other hand, for the remaining centers $x_i^n$, 
we  transport the particles from $\P_n$ to $\P_{n-1}$, using the branches of the dyadic irrigation plan $\chi_n$, defined as in Lemma \ref{l:42}. Notice on each such branch, $m_i^n < z_0$.  Then from (\ref{nw2}) and (\ref{1020}), the weight $W^n_{n,i}:[0,{\sqrt{d}L\over 2^{n+1}}]\mapsto \R_+$ on the branch $\gamma_i$ from $\P_n$ to $\P_{n-1}$ satisfies
\bel{1023}  W^n_{n,i}(s)~=~\Big((m_i^n)^{1-\beta} + c(1-\beta){\sqrt{d}L\over 2^{n+1}} \Big)^{1\over 1-\beta}~\leq~(z_0^{1-\beta} + z_0^{1-\beta})^{1\over 1-\beta}~=~2^{1\over 1-\beta}z_0,  \eeq
where we compute the weight $W_{n,i}^n$ as solution to $\dot W_{n,i}^n = c(W_{n,i}^n)^\beta$. Let $W_i$ be the corresponding solution of (\ref{qs1}) with $m_i(t)$ replaced by constant multiplicity $m_i^n$, by (\ref{fp5}) and comparision principle from ODE theory, 
\bel{121a}W_i(s)~\leq~W_{n,i}^n(s). \eeq
 
Then clearly the total cost on these dyadic branches from $\P_n$ to $\P_{n-1}$ is bounded by $E_n^n$, given in (\ref{e1}).
\v
{\bf 2.} From $\P_{n-1}$ to $\P_{n-2}$. After removing the point masses transported by branches in  $\S_{n}$, we still denote the remaining measure as $\mu_n$, and transport $\mu_n$ to the centers in $\P_{n-1}$, using the branches from $\P_n$ to $\P_{n-1}$ of the dyadic irrigation plan $\chi_n$. Notice that after removing the masses transported by branches in $\S_n$, $m_i^n\leq z_0$ for each $1\leq i\leq 2^{nd}$, with some $m_i^n=0$.

 For each center $x_i^{n-1}$ in $\P_{n-1}$, when 
\bel{1024} \sum_{j\in \O(i)} m_j^n ~\geq~z_0 \eeq
we then connect $x_i^{n-1}$ to the origin directly by a straight branch. Let $\S_{n-1}$ be the set of all such branches. Similarly as in (\ref{1021}),  the weight on each branch in $\S_{n-1}$ is bounded by $e^{c\sqrt{d}L}M$. For the remaining $x_i^{n-1}$, we transport the flux from $\P_{n-1}$ to  $\P_{n-2}$, by the branches of dyadic irrigation plan $\chi_n$. From (\ref{nw13}) and (\ref{fp5})-(\ref{1020}), on each dyadic branch $\gamma_i$ from $\P_{n-1}$ to $\P_{n-2}$,
\bel{1025} W_i(s)~\leq~W_{n-1,i}^n(s)~<~2^{1\over 1-\beta}z_0,\qquad\qquad s\in [0,\, {\sqrt{d}L\over 2^n}]. 
 \eeq
 Then clearly the total cost on these dyadic branches from $\P_{n-1}$ to $\P_{n-2}$ is bounded by $E_{n-1}^n$, defined by (\ref{e2}).
\v
{\bf 3.} By backward induction we construct the irrigation plan until to the level $\P_{n_0}$. For each $k> n_0$, from $\P_{k}$ to $\P_{k-1}$, there are two types of paths, one is the branches in $\S_k$, and the other one is the dyadic branches of $\chi_n$.  Clearly we have
\bel{1026} \#\Big( \bigcup_{k>n_0}^n \S_k \Big)~\leq~{M\over z_0}  \eeq
where $M$ is the total mass of $\mu$. Indeed, from our construction, each branch in $\cup_{k>n_0}^n \S_k$ will transport distinct groups of particles with mass $\geq z_0$, the total mass of $\mu_n$ is $M$, thus we have the upper bound in (\ref{1026}).
For each branch in $\S_k$, there is an uniform bound (\ref{1021}) on the weight $W(t)$, and the length of each branch is bounded by $\sqrt{d}L$, thus the total cost $J$ on branches in  $\S_k, k> n_0 $ is bounded by 
\bel{1027} J~\leq~{M\over z_0} \cdot \Big(e^{c\sqrt{d}L M}\Big)^\alpha \sqrt{d}L ~\doteq~\kappa_0  \eeq 
On the other hand, the total cost  $I$ on the dyadic branches is bounded by
\bel{1028}  I~\leq~\sum_{k>n_0}^n E_k^n~\leq~C\Big( M^\alpha L + L^{1 + {\alpha\over 1 - \beta}}  \Big)~\doteq~\kappa_1 \eeq	
where the last inequality comes from (\ref{univ1}).

Notice the bounds in (\ref{1027})-(\ref{1028}) are independent of $n$, therefore, there exists a uniform constant $C>0$, such that for each dyadic approximation $\mu_n$, we have $\I^{W,\alpha}(\mu_n)\leq C$. Thanks to  Theorem \ref{t:meas}, we conclude that $\I^{W,\alpha}(\mu)~\leq~C$.
\endproof

\subsection{Examples of non-irrigable measures}
In the following we show some cases for measures $\mu$ 
with infinite weighted irrigation cost $\I^{W,\alpha}$. 

\begin{definition}\label{d:alf}
	Let $\mu$ be a positive, bounded measure in $\R^d$. If there exists $\gamma >0$ and a constant $C\geq 1$ such that
	\bel{1001} {1\over C} r^\gamma~\leq~\mu(B(x,r))~\leq~Cr^\gamma,\qquad \forall x\in \hbox{supp}(\mu),~r\in [0,1], \eeq 
   then we say $\mu$ is Ahlfors regular in dimension $\gamma$. Here $\hbox{supp}(\mu)$ is the support of $\mu$, $B(x,r)$ is the ball of radius $r$ that centered at $x$. 
	\end{definition}

\begin{remark}
	If a measure $\mu$ is Ahlfors regular in dimension $\gamma$, then one can prove $supp(\mu)$ has Hausdorff dimension $\gamma$. Indeed, consider any covering $\cup_{i=1}^\infty B(x_i,r_i)$ of $supp(\mu)$, consists of closed balls with radius less than 1.
	 From the second inequality in (\ref{1001}) one has
	$$\sum_{i\geq 1}({r_i })^\gamma~\geq~\sum_{i\geq 1}{\mu(B(x_i,r_i))\over C}~\geq~{M\over C}~>~0,$$
	which implies $supp(\mu)$ has Hausdorff dimension $\geq\gamma$. On the other hand, by the Vitali's Convering Theorem\cite{EVANS}, there exists a  countable subcollection of disjoint $B(x_i,r_i)$, which we still denote as $\sum_{i=1}^\infty B(x_i,r_i)$, such that  $supp(\mu)\subseteq \cup_{i=1}^\infty B(x_i,5r_i)$. Then from the first inequality in (\ref{1001}), since $B(x_i,r_i)$ are disjoint,
	$$\sum_{i\geq 1}(5 r_i)^\gamma~=~5^\gamma\sum_{i\geq 1} r_i^\gamma~\leq~5^\gamma C \sum_{i\geq 1}\mu(B(x_i,r_i))~\leq~5^\beta C M,$$
	and it implies the Hausdorff dimension of $supp(\mu)$ $\leq \gamma$.
\end{remark}

For the irrigation cost $\I^\alpha(\cdot)$ without weights that defined in \cite{MMS}, we recall the following theorem. For a proof, 
see  Theorem 1.2 in \cite{MMS}.

\begin{theorem}\label{t:old}
	
	Let $\mu$ be a finite $\alpha$-irrigable measure, with $\alpha\in (0,1)$. That is, $\I^\alpha(\mu)<\infty$. Then there is a Borel set $E\subseteq \R^d ,\mu(\R^d\setminus E)~=~0$,  such that for any $s> {1\over 1 - \alpha}$,
	$$ \H^s(E)~=~0,   $$
	where $\H^s(E)$ is the $s$-Hausdorff measure of the set $E$.
In other words, if $\mu$ is $\alpha$-irrigable, then $\mu$ is concentrated on a  set $E$ with Hausdorff dimension $\leq {1\over 1-\alpha}$.
\end{theorem}

\begin{remark}\label{r:hd}
	As mentioned in \cite{BSUN2}, for any bounded Radon measure $\mu$, we always have $\I^{W,\alpha}(\mu)\geq \I^{\alpha}(\mu).$ Therefore, if $\I^{W,\alpha}(\mu)<+\infty$, from Theorem \ref{t:old}, $\mu$ is concentrated on a set $E$ with Hausdorff dimension $\leq {1\over 1-\alpha}$.
\end{remark}

\begin{lemma}\label{l:mout}
	
	If $\mu$ is a bounded Radon measure as in Theorem \ref{t:i} and let $\chi$ be an irrigation plan of $\mu$ with finite weighted irrigation cost $\E^{W,\alpha}(\chi)<\infty$. Then for any $r>0$, 
	\bel{1004} \mu(\R^d\setminus B(0,r) )~\leq~\left({\E^{W,\alpha}(\chi)\over r} \right)^{{1\over \alpha}}.  \eeq
\end{lemma}
\v
{\bf Proof.} The function 
$$ x~\mapsto~ \left(x^{1-\beta} + c(1 - \beta)(r - t) \right)^{{1\over 1 - \beta}},\qquad x\in \R_+$$
is concave. Let $m_r\doteq \mu(\R^d\setminus B(0,r))$, then by definition (\ref{10220}) and (\ref{qs9})  we have
\bel{1006} \int_0^r \left( m_r^{1-\beta} + c(1-\beta)(r-t) \right)^{{\alpha\over 1 -\beta}}\, dt~\leq~\E^{W,\alpha}(\chi).      \eeq

Since $r-t\geq 0$, it implies that
$$ \left(m_r^{1-\beta}\right)^{\alpha\over 1-\beta}\cdot r~=~m_r^\alpha\cdot r~\leq~\E^{W,\alpha}(\chi), $$
which completes the proof of (\ref{1004}).


\v

\begin{theorem}\label{t:ni}
	Let $\mu$ be a positive, bounded Radon measure in $\R^d$  and  Ahlfors regular in dimension $d$.   
	Let $f$ satisfy {\bf (A1)}.   
	
	If either $\alpha<1-{1\over d}$ or 
	\bel{fp7}
	\liminf_{z\to 0^+} ~z^{-\beta} f(z)~>~0\eeq
	for some $\beta<1-{1\over d-1}$, then $\I^{W,\alpha}(\mu)=+\infty$.
\end{theorem}

{\bf Proof.}  

CASE 1: If $\alpha < 1 - {1\over d}$, then ${1\over 1-\alpha}< d$.  Suppose $\I^{W,\alpha}(\mu)<+\infty$, by Remark \ref{r:hd}, $\mu$ is concentrated on a set $E$ with Hausdorff dimension $\leq {1\over 1-\alpha}<d$, which is a contradiction to the assumption that $\mu$ is Ahlfors regular in dimension $d$ (see Remark 3.5). Thus, we have $\I^{W,\alpha}(\mu)=+\infty$.
\v
CASE 2: The assumption (\ref{fp7}) implies that, for some constants $c,z_0 >0$,
\bel{fp8} f(z)~\geq~c z^\beta\qquad\qquad \forall z\in [0, z_0].\eeq
   Since $\mu$ is Ahlfors regular in dimension $d$, then for each irrigation plan $\chi$ and any $\delta>0$, there are $\O({1\over \delta^d})$ disjoint cubes with diameter $\delta$ and each of them has measure $\approx \delta^d$. In each cube, the lower bound for the cost is
\bel{1009} \int_0^\delta \left( \delta^{d(1-\beta)} + c(1 - \beta)(\delta - t)  \right)^{{\alpha\over 1 - \beta}}\, dt   \eeq 
and the total number of such disjoint cubes is ${1\over \delta^d}$. 
\bel{1010}\bega{rl}\ds \E^{W,\alpha}(\chi)&\ds\geq~{1\over \delta^d}  \int_0^\delta \left( \delta^{d(1-\beta)} + c(1 - \beta)(\delta - t)  \right)^{{\alpha\over 1 - \beta}}\, dt\\[3mm]
&\ds \geq~{1\over \delta^d} \int_0^\delta \Big( c(1-\beta)(\delta - t) \Big)^{\alpha\over 1-\beta}\, dt\\[3mm]
 &\ds \geq~C\delta^{{1+\alpha-\beta\over 1-\beta}-d} \enda \eeq
 where $C$ is some constant independent of $\delta$. Since $1\geq\alpha> 1-{1\over d}, 1-{1\over d-1}>\beta >0$, we have ${1+\alpha-\beta\over 1-\beta}<d$.
  Sending $\delta$ to $0+$, the right hand side in (\ref{1010}) goes to $+\infty$. Thus, for any irrigation plan $\chi$ of $\mu$, $\E^{W,\alpha}(\chi)=+\infty$, thus $\I^{W,\alpha}(\mu)=+\infty$.

\v

{\bf Acknowledgement   } This research was partially supported by NSF grant DMS-1714237, ``Models of controlled biological growth".  The author wants to thank his thesis advisor   Professor Alberto Bressan for his many useful comments and suggestions. The author also thanks the anonymous referees for their time and helpful suggestions.

\end{document}